\documentclass[12pt]{article}

\usepackage{amsmath}
\usepackage{amssymb}

\newcommand{\RR}{\mathbb{R}}
\newtheorem{theorem}{Theorem}[section]

\newtheorem{corollary}{Corollary}[section]
\newtheorem{example}{Example}[section]

\newtheorem{definition}{Definition}[section]

\begin{document}

\title{On Lyapunov-Krasovskii Functionals for Switched Nonlinear 
Systems with Delay}

\author{A. Yu. Aleksandrov \thanks{Faculty of Applied Mathematics and Control Processes, St. Petersburg State University,  St. Petersburg, Russia; email: alex43102006@yandex.ru; supported by by the St. Petersburg State University, project 
no. 9.38.674.2013, by the Russian Foundation of Basic Researches, grants no. 13-01-00376-a and 13-08-00948-a}
\and Oliver Mason\thanks{Hamilton Institute, National University of Ireland Maynooth, Co. Kildare, Ireland; email: oliver.mason@nuim.ie; Supported by the Irish Higher Educational Authority, PRTLI 4 Network Mathematics Grant}}

\maketitle

\begin{abstract}
We present a set of results concerning the existence of Lyapunov-Krasovskii functionals for classes of nonlinear switched systems with time-delay.  In particular, we first present a result for positive systems that relaxes conditions recently described in \cite{SunWang} for the existence of L-K functionals.  We also provide related conditions for positive coupled differential-difference positive systems and for systems of neutral type that are not necessarily positive.  Finally, corresponding results for discrete-time systems are described.  
\end{abstract}


\bigskip

\section{Introduction}
\label{sec:Intro}
In the stability analysis of positive systems (\cite{FR00}, \cite{Kac}), 
\emph{linear} Lyapunov functions have received considerable attention recently (\cite{AP}, \cite{KMS09}, \cite{ELENA1}, \cite{ELENA2}).  Initial efforts in this direction considered linear time-varying systems and switched linear systems in particular.  More recently, researchers have focussed on deriving usable stability conditions for classes of nonlinear systems and systems subject to time-delay.  Our results are concerned with the use of linear Lyapunov-Krasovski (L-K) functionals for classes of nonlinear time-delay systems.  In this introduction, 
we first recall some relevant work in order to place our results in context.

The notation of the paper is very standard.  We denote the space of continuous functions on $([-\tau, 0]$ taking values in $\mathbb{R}^n$ by $C([-\tau, 0], \mathbb{R}^n)$.  A matrix $A \in \mathbb{R}^{n \times n}$ is Hurwitz if all of its eigenvalues have negative real part.  $A$ is said to be Schur-Cohn if all of its eigenvalues have modules less than one.  We use $e$ to denote the vector all of whose entries are equal to one, $e = (1, \ldots ,  1)^T$.  For a vector $v \in \mathbb{R}^n$, $v \geq 0$ means that $v_i \geq 0$ for $1 \leq i \leq n$.  Similarly, $v > 0$ means that $v \geq 0$, $v \neq 0$, while $v \gg 0$ means that $v_i > 0$ for all $i$.   For a vector $x \in \mathbb{R}^n$, the notation $|x|$ denotes the vector $(|x_1|, \ldots, |x_n|)^T$.  For a functional $V$ defined on $C([-\tau, 0], \mathbb{R}^n)$, $\partial^+ V$ denotes the upper-right Dini derivative of $V$ \cite{Hale}.

It is important to note that the interest in using linear Lyapunov functions to study positive systems stretches back beyond the last decade.  In fact, work of a similar spirit can be found in \cite{Persid}.  In this paper, the system
\begin{equation}\label{eq:Pers1}
\dot {x}(t)=A f(x(t))
\end{equation}
is considered.  The nonlinearity $f:\mathbb{R}^n \rightarrow \mathbb{R}^n$, is continuous, diagonal 
$$f(x_1, \ldots, x_n) = (f_1(x_1), f_2(x_2), \ldots, f_n(x_n))^T$$
and assumed to satisfy
\begin{equation}\label{eq:nonlin}
x_if_i(x_i)>0 \ \ \ \text{for}\ \ x_i\neq 0. 
\end{equation}
Furthermore, the matrix $A$ is assumed to be \emph{Metzler}, meaning that its 
off-diagonal elements $a_{ij}$, $i \neq j$ are nonnegative.  
Taken together, these conditions mean that (\ref{eq:Pers1}) is a positive system with an equilibrium at the origin.  It is established in \cite{Persid} that if $A$ is a Hurwitz matrix then the zero solution of (\ref{eq:Pers1})  is asymptotically stable for any admissible nonlinearity $f(x)$, and a Lyapunov function can be
chosen in the form 
$$
V(x)= \nu^T |x| = \sum_{i=1}^n \nu_i |x_i|
$$
where the positive vector $\nu=(\nu_1,\ldots,\nu_n)^T$  is a solution of 
the system   $A^T\nu \ll 0.$

More recently, in \cite{Haddad} linear L-K functionals were used to analyse 
positive linear time-delay
systems of the form
$$
\dot {x}(t)= A x(t)+  B x(t -\tau).
$$
This system is positive if and only if $A$ is Metzler and $B$ is nonnegative.  
The main result of \cite{Haddad} established that the zero solution of this 
system is asymptotically stable for all $\tau \geq 0$ if and only if there is 
some vector $\nu$ with positive entries such that $(A+B)^T \nu$ is entrywise negative.  Moreover, it was shown that in this case, a L-K functional of the form
$$
V(x_t)= \nu^T |x(t)|+ 
\int_{t-\tau}^t  \nu^T  Bx(\theta) d\theta
$$                                                                    
exists.  As usual $x_t$ denotes the continuous function on $[-\tau, 0]$ given 
by $x_t(\theta) = x(t + \theta)$ for $\theta \in [-\tau, 0]$.  A simple 
nonlinear extension of this result, where the nonlinearity occurs only in 
the delay term and is assumed to be sublinear is also presented.  Where there is no risk of confusion, we shall refer to the stability of the system rather than of the zero solution.  

An extension of the core result of \cite{Persid} to switched systems was 
given in \cite{AP} (related results can be found in \cite{ACPZ11}).  This paper considers a switched system
$\dot {x}(t)=A^{(\sigma)} f(x(t))$ with a corresponding family of constituent systems
$\dot {x}(t)=A^{(s)} f(x(t)), \quad s=1,\ldots,N$
where $A^{(s)}$  are Metzler and Hurwitz matrices, and
the diagonal nonlinearity $f:\RR^n \rightarrow \RR^n$, 
satisfies the condition (\ref{eq:nonlin}).  Extending the result of 
\cite{Persid}, it is shown that if there exists a positive vector $\nu \gg 0$ 
satisfying ${A^{(s)}}^T\nu  \ll 0, \quad s=1,\ldots,N $,
then the nonlinear switched system is asymptotically
stable for any admissible nonlinearity $f(x)$ and for arbitrary switching signals. Moreover, 
a common Lyapunov function can be 
chosen in the form $\nu^T |x|$.

Switched nonlinear systems with time-delay 
\begin{equation}\label{eq:DelSys1}
\dot {x}(t)=A^{(\sigma)} f(x(t))+B^{(\sigma)} f(x(t-\tau))
\end{equation} 
were considered in \cite{SunWang}.  The corresponding family of subsystems is given by
\begin{equation}
\label{eq:Fam1}
\dot {x}(t)=A^{(s)} f(x(t))+B^{(s)} f(x(t-\tau)), \quad s=1,\ldots,N.
\end{equation}
Here $A^{(s)}$  are Metzler and Hurwitz matrices,
$B^{(s)}$  are nonnegative matrices, and
the diagonal nonlinearity $f:\RR^n \rightarrow \RR^n$, 
is continuous and satisfies (\ref{eq:nonlin}).

Throughout the paper, an \emph{admissible switching signal} is a piecewise constant function $\sigma:[0, \infty) \rightarrow \{1, \ldots, N\}$, which is right continuous and has only finitely many discontinuities on every bounded interval.  

We say that (\ref{eq:DelSys1}) is \emph{absolutely stable} if the origin is a globally asymptotically stable equilibrium for every admissible switching signal $\sigma$, every continuous diagonal nonlinearity satisfying (\ref{eq:nonlin}) and every nonnegative delay $\tau$.  For a non-switched system, absolute stability is defined with respect to the family of nonlinearities and all nonnegative delays only. 

Let $b^{(s)}_{ij}$ denote the entries of $B^{(s)}$, and write $\bar B=\{\bar b_{ij}\}$, where
$$
\bar b_{ij}=\max_{s=1,\ldots,N} b^{(s)}_{ij}.
$$
In \cite{SunWang}, it was proved that if there exists a positive vector $\nu=(\nu_1,\ldots,\nu_n)^T$ such that
$\left({A^{(s)}}+\bar B\right)^T\nu\ll 0$ for $s=1,\ldots,N$,
then the switched system (\ref{eq:DelSys1}) is 
absolutely stable,
and a common L-K functional can be
chosen in the form 
$$
V(x_t)=\sum_{i=1}^n \nu_i  |x_i(t)|+ 
\sum_{i=1}^n \nu_i \int_{t-\tau}^t \sum_{j=1}^n  \bar b_{ij} |f_j(x_j(\theta))|
d\theta. 
$$                                                                    
Our goal is to obtain less restrictive stability conditions for 
(\ref{eq:DelSys1})
than those in \cite{SunWang}.
Moreover, we will extend our approach to other distinct classes of 
switched systems including systems with several delays, coupled differential and
difference systems, and neutral type systems.

\section{Switched Nonlinear Differential Systems with Time-Delay}
\label{sec:swit}
In this section, we describe sufficient conditions for the absolute stability of the system (\ref{eq:DelSys1}) using a form of L-K functional which differs from those used in (\cite{Haddad}, \cite{SunWang}).
To begin, we consider the single delay case in the
interest of notational simplicity and clarity; the arguments for the several delay case are identical.

\subsection{A Single Delay}

\begin{theorem}\label{thm:Swit1}
Consider the system (\ref{eq:DelSys1}).  Assume that $A^{(s)}$ is Metzler and $B^{(s)}$ is nonnegative for $s = 1, \ldots, N$.  Assume that there exists  
a vector $\nu \gg 0$ satisfying the
inequalities
\begin{equation}\label{eq:swit1}
(A^{(s)}+B^{(r)})^T \nu \ll 0, \qquad s,r=1,\ldots ,N.
\end{equation}
Then, the switched system (\ref{eq:DelSys1}) is absolutely stable; in fact, there exist positive constants $\mu_1, \ldots, \mu_n$ such that 
\begin{equation}
\label{eq:LK1}
V=\sum_{i=1}^n \nu_i |x_i(t)|+\sum_{i=1}^n \mu_i \int_{t-\tau}^t 
|f_i(x_i(z))|dz
\end{equation}
defines a L-K functional for the family (\ref{eq:DelSys1}).
\end{theorem}
                    
\textbf{Proof:}  
By assumption, there exists some vector $\nu \gg 0$ with $(A^{(s)} + B^{(r)})^T \nu \ll 0$ for all $s, r$.  Choose some positive vector $q \gg 0$ such that
\begin{equation}
\label{eq:Swit11} (A^{(s)} + B^{(r)})^T \nu \leq - q.
\end{equation} 
Next set $w = \textrm{max}(B^{(s)})^T \nu$, where the maximum is taken elementwise.  Clearly $w \geq 0$.  Now choose $\mu = w + \frac{q}{2}$.  It follows readily from (\ref{eq:Swit11}) that for all $s \in \{1, \ldots , N\}$, 
\begin{equation}
\nonumber (A^{(s)})^T \nu + w \leq  - q
\end{equation}
so that 
\begin{equation}
\label{eq:Swit12} 
(A^{(s)})^T \nu \leq - q - w \ll - \mu.
\end{equation}
On the other hand, it is immediate from the definition of $w$ that for all $s \in \{1, \ldots  , N\}$, 
\begin{equation}
\label{eq:Swit13} (B^{(s)})^T \nu \ll \mu.
\end{equation}
If we differentiate (\ref{eq:LK1}) along any of the systems (\ref{eq:Fam1}) we find that the derivative satisfies
\begin{eqnarray*}
\label{eq:Swit14} 
\partial^+ V  \leq  |f(x(t))|^T({A^{(s)}}^T \nu + \mu) 
\\
+ 
|f(x(t-\tau))|^T\left({B^{(s)}}^T \nu - \mu\right)  
\leq-\beta \sum_{j=1}^n |f_j(x_j(t))|,
\end{eqnarray*}
where $\beta$ is a positive constant.
Thus $V$ defines a L-K functional for (\ref{eq:DelSys1}) as claimed.  As this is true for any choice of $f$ and $\tau$, it follows that (\ref{eq:DelSys1}) is absolutely stable as claimed.

\textbf{Remark:} The results of \cite{SunWang} can be applied to (\ref{eq:DelSys1}) to provide a closely related set of conditions for absolute stability.  However, we use a slightly different type of L-K functional and moreover our conditions are less restrictive than those in \cite{SunWang}.  The following example illustrates this last point.

\begin{example}
\label{ex:swit2}  Consider an example of the system (\ref{eq:DelSys1}) 
 with $n=2$, $N = 2$ and 
$$A=A_1 =A_2= \left(\begin{array}{c c}
				-2 & 0 \\
				0 & -2 
				\end{array}\right),
$$
$$				
B_1 = \left(\begin{array}{c c}
				1 & 1 \\
				1 & 0 
				\end{array}\right), \qquad
B_2 = \left(\begin{array}{c c}
				0 & 1 \\
				3 & 0 
				\end{array}\right).
$$				

On the one hand, if we apply the approach from \cite{SunWang}, we would first construct the matrix
$$\bar B= \left(\begin{array}{c c}
				1 & 1 \\
				3 & 0 
				\end{array}\right),
$$
which is the entrywise maximum of $B_1$ and $B_2$.  In this case, it can be easily verified that there is no $\nu \gg 0$ such that $(A+\bar B)^T \nu\ll 0.$  Hence, we cannot apply Theorem 2 of \cite{SunWang} in this case.

On the other hand, applying Theorem \ref{thm:Swit1}, we can check that
$$
(A+ B^{(s)})^T \nu \ll 0, \quad s=1,2, 
$$
with $\nu = (7, 4)^T$ for example.  
\end{example}

\textbf{Remark:}  It is worth mentioning that the existence
of a common Lyapunov function in the form 
$$
V(x)= \nu^T |x|
$$
for the corresponding delay free system does  not imply
the existence of a common L-K functional in the form (\ref{eq:LK1})
for the system (\ref{eq:DelSys1}).
              
\begin{example}
\label{ex:swit3}  Consider the simple 1-dimensional systems
$$
\dot x(t) = -0.2 x(t) + 0.1 x(t-\tau),  
$$
$$
\dot x(t) = -0.9 x(t) + 0.8 x(t-\tau).
$$  
As the undelayed systems obtained by setting the delay to zero are 
identical, they clearly have a common Lyapunov function of the required form.  
However,  it can be easily verified that no common L-K functional in the form
(\ref{eq:LK1})  exists for the switched systems.
\end{example}

\subsection{Several Delays}

Finally, for this section, consider the switched system with several delays
\begin{equation}\label{eq:MulDel}
\dot x(t)=        A^{(\sigma)}f(x(t))+B^{(\sigma)}_1f(x(t-\tau_1))
\end{equation}
$$
 +\ldots+
B^{(\sigma)}_lf(x(t-\tau_l))        
  $$
and the corresponding family of subsystems
\begin{equation}\label{eq:MulDelFam}
\dot x(t)=        A^{(s)} f(x(t))+B^{(s)}_1 f(x(t-\tau_1))+\ldots
\end{equation}
$$
+
B^{(s)}_lf(x(t-\tau_l)), \ \ \ s=1,\ldots,N,
$$    
  where $A^{(s)}$ are Metzler matrices, and $B^{(s)}_1,\ldots,B^{(s)}_l$ 
are nonnegative matrices for $1 \leq s \leq N$,  
and the nonlinearities $f(x)=(f_1(x_1),\ldots,f_n(x_n))^T$ are continuous 
and satisfy the 
condition (2).  The real numbers $\tau_1,\ldots,\tau_l$ are nonnegative delays.

It is relatively straightforward to adapt the argument 
used in the proof of Theorem 2.1 to obtain the following result 
for systems with several delays.

\begin{theorem}\label{thm:mult}
Assume that there exists a positive vector $\nu\gg 0$ satisfying the
inequalities
\begin{equation}\label{eq:mult1}
(A^{(s)}+B^{(p_1)}_1+\ldots+B^{(p_l)}_l)^T \nu \ll 0,  
\end{equation}
for $ s,p_1,\ldots,p_l \in \{1,\ldots,N\}.$
Then the switched system (\ref{eq:MulDel}) is absolutely stable, and for the family (\ref{eq:MulDelFam})
there exist positive constants $\mu_{rj}$ for $r = 1, \ldots, l$, $j = 1, \ldots, n$ such that  
\begin{equation}\label{eq:LKMul}
V=\sum_{i=1}^n \nu_i |x_i(t)|+\sum_{i=1}^n  \sum_{r=1}^l \int_{t-\tau_r}^t 
\mu_{ri}|f_i(x_i(z))|dz
\end{equation}
is a common L-K  functional for the system (\ref{eq:MulDel}).
\end{theorem}

\textbf{Remark:} To check for the existence of a positive vector $\nu\gg 0$ 
satisfying the inequalities  (\ref{eq:mult1}) and to find this vector, one can apply approaches proposed in 
(\cite{KMS09}, \cite{ACPZ11}, \cite{Siljak}).

\section{Coupled Nonlinear Differential and Difference Systems }
\label{sec:coup}
Systems described by coupled differential and difference equations arise in applications hydraulics, circuit theory, and the analysis of partial differential equations.  They can also be used to reduce the dimensionality of a large-scale delayed system, by suitably exploiting the fact that, while an overall systems may be high-dimensional, delays often occur in scalar or low-dimensional parts of the system \cite{GuLiu}.  In this section, we present sufficient conditions for absolute stability of a class of nonlinear coupled differential-difference systems using a L-K functional of linear form.

\subsection{Nonswitched Systems}

Formally, we consider the system
\begin{eqnarray}
\label{eq:Coupled1}
\dot{x}(t) &=& A f(x(t)) + By(t - \tau),\\
\label{eq:Coupled2}
y(t)&=& Cf(x(t)) + Dy(t-\tau).
\end{eqnarray}
Here $x(t)\in \mathbb{R}^n$, $y(t) \in \mathbb{R}^m$, $A \in \mathbb{R}^{n \times n}$, $B \in \mathbb{R}^{n \times m}$, $C\in \mathbb{R}^{m \times n}$, $D \in \mathbb{R}^{m \times m}$.
The vector field $f:\mathbb{R}^n \rightarrow \mathbb{R}^n$ is a continuous 
diagonal nonlinearity $f(x)=(f_1(x_1),\ldots,f_n(x_n))^T$ satisfying (\ref{eq:nonlin}).  Such a vector field is said to be an \emph{admissible nonlinearity}.  

Initial conditions for (\ref{eq:Coupled1}), (\ref{eq:Coupled2}) are specified by a vector $x(0) \in \mathbb{R}^n$ together with a function $y_0$  in $C([-\tau, 0], \mathbb{R}^n)$.  We write $\|x(0), y_0\|$ for $\textrm{max}(\|x(0)\|_2, \|y_0\|_{\infty})$ where $\|\cdot \|_2$ is the usual Euclidean norm on $\mathbb{R}^n$ and $\|\cdot \|$ is the $l_{\infty}$ norm on $C([-\tau, 0], \mathbb{R}^n)$.  Under our assumptions on the vector field $f$, 
the origin is always an equilibrium point of (\ref{eq:Coupled1}), (\ref{eq:Coupled2}).  
We recall the following stability definitions from \cite{GuLiu}.  

\begin{definition}
\label{def:stable}
The origin is a stable equilibrium of (\ref{eq:Coupled1}), (\ref{eq:Coupled2}) if for any $\epsilon > 0$, there exists some $\delta > 0$ such that $\|x(0), y_0\| < \delta$ implies that $\|x(t), y_t\| < \epsilon$ for all $ t \geq 0$.  It is said to be asymptotically stable if  it is stable and, in addition, $\|x(t), y_t\| \rightarrow 0$ as $t \rightarrow \infty$ for all initial conditions $x(0), y_0$.  
\end{definition}

Our primary concern throughout is with asymptotic stability.  In an abuse of 
notation, we shall refer to the stability of the system 
(\ref{eq:Coupled1}), (\ref{eq:Coupled2}) rather than the stability of the equilibrium at the origin. 

The system (\ref{eq:Coupled1}), (\ref{eq:Coupled2}) is said to be \emph{absolutely stable} if its zero solution
is asymptotically stable for every admissible nonlinearity and  every $\tau\geq 0$.

{\bf Remark:} As noted in \cite{GuLiu}, the second equation (\ref{eq:Coupled2}) can be viewed as a separate system in which $y_t$ is the state and $x(t)$ acts as an input.  Moreover, it is clear that for the overall system (\ref{eq:Coupled1}), (\ref{eq:Coupled2}) to be asymptotically stable, the system described by (\ref{eq:Coupled2}) must be input-to-state stable.  This immediately  implies that Schur-stability of the matrix $D$ is a necessary condition for asymptotic stability of (\ref{eq:Coupled1}), (\ref{eq:Coupled2}).

Furthermore, if (\ref{eq:Coupled1}), (\ref{eq:Coupled2}) is absolutely stable,  then it is asymptotically stable 
when $f_i(x_i)=x_i$ for all $i$ and $\tau=0$.  It follows immediately that the matrix
\begin{equation}
\label{eq:Coupled3}
 A+B(I-D)^{-1}C 
\end{equation}
must be Hurwitz.
                 

Our next result considers a \emph{positive} coupled differential-difference system and shows that in this case, the last condition above is sufficient for the existence of a linear L-K functional.

\begin{theorem}
\label{thm:Coupled1}
Consider the system (\ref{eq:Coupled1}), (\ref{eq:Coupled2}).  Assume that $A$ is a Metzler matrix and that 
$B, C, D$ are nonnegative matrices.  Furthermore, assume that  
$D$ is Schur-Cohn stable and that (\ref{eq:Coupled3})
is  Hurwitz stable.
Then there exist positive constants $\mu_i$, $\nu_i$, $1 \leq i \leq n$, such that
\begin{equation}
\label{eq:LKCoupled}
V=\sum_{i=1}^n \nu_i |x_i(t)|+ \sum_{i=1}^n \mu_i
\int_{t-\tau}^t |y_i(z)|dz
\end{equation}
is a L-K functional for the system defined by (\ref{eq:Coupled1}), (\ref{eq:Coupled2}).
\end{theorem}

\textbf{Proof:} 
As $D$ is a Schur-Cohn matrix, there exists some vector $v \gg 0$ with $(D - I)^T v \ll 0$.  Furthermore, it follows that $I-D$ is an M-matrix and hence has a non-negative inverse.  Therefore, the matrix $B(I-D)^{-1}C$ is also non-negative and $A + B(I-D)^{-1}C$ is Metzler.  By assumption, this last matrix is also Hurwitz and thus we can choose some $\nu \gg 0$ such that 
\begin{equation}\label{eq:nu}
(A+B(I-D)^{-1}C)^T\nu\ll 0.
\end{equation}
Next set 
\begin{equation}\label{eq:mu}
\mu=(I-D^T)^{-1} B^T \nu+\varepsilon v
\end{equation}
where $\varepsilon$ is a positive constant.
It follows readily from (\ref{eq:mu}) that 
\begin{equation}\label{eq:B}
B^T \nu + (D-I)^T \mu = \varepsilon (D-I)^T v \ll 0
\end{equation}
for any $\varepsilon > 0$.  It also follows from (\ref{eq:mu}) that
\begin{equation}\label{eq:A1}
A^T \nu + C^T \mu =  (A+B(I-D)^{-1}C)^T\nu + \varepsilon C^T v .
\end{equation}
From (\ref{eq:nu}), we see that we can choose $\varepsilon > 0$ small enough to ensure that 
\begin{equation}\label{eq:A}
A^T \nu + C^T \mu \ll 0.
\end{equation}

Next note that the Dini upper-right derivative of the functional (\ref{eq:LKCoupled}) satisfies
the inequality
\begin{equation}\label{eq:DiniCoupled}
\partial^+ V\leq \sum_{j=1}^n   |f_j(x_j(t))| \left( \sum_{i=1}^n  \nu_i a_{ij} +
\sum_{i=1}^n  \mu_i c_{ij}\right)  
\end{equation}
$$
+ \sum_{j=1}^n   |y_j(t-\tau)| \left( \sum_{i=1}^n  \nu_i b_{ij} +
\sum_{i=1}^n  \mu_i d_{ij}-\mu_j\right).
$$
We can write the right hand side of (\ref{eq:DiniCoupled}) as 
\begin{equation}\label{eq:DiniCoupled2}
|f(x)|^T (A^T \nu + C^T \mu) + |y(t - \tau)|^T( B^T \nu + (D-I)^T \mu).
\end{equation}
From (\ref{eq:A}), it follows that we can find some $\beta  > 0$ such that $A^T\nu + C^T \mu \leq -\beta e$ where $e = (1, 1, \ldots, 1)^T$.  Combining this with (\ref{eq:B}) and (\ref{eq:DiniCoupled2}) yields
\begin{equation}
\partial^+ V\leq -\beta \sum_{j=1}^n   |f_j(x_j(t))|.
\end{equation}
It now follows immediately from Theorem 3 of \cite{GuLiu} that 
(\ref{eq:LKCoupled}) is a L-K functional for the system (\ref{eq:Coupled1}), (\ref{eq:Coupled2}) and that the system is absolutely stable as claimed. 

If we combine Theorem \ref{thm:Coupled1} with the remark preceding it, we obtain the following result. 
\begin{corollary}\label{cor:Coupled1}
Consider the system described by (\ref{eq:Coupled1}), (\ref{eq:Coupled2}).  Assume that $A$ is a Metzler matrix,  
 $B, C, D$ are non-negative matrices, and 
$D$ is Schur-Cohn stable.
The system (\ref{eq:Coupled1}), (\ref{eq:Coupled2}) is absolutely stable if and only if the matrix (\ref{eq:Coupled3})
is Hurwitz stable.
\end{corollary}

\subsection{Switched Systems}

We next consider the switched coupled system
\begin{eqnarray}\label{eq:SwCoupled}
\dot{x} &=& A^{(\sigma)}  f(x(t)) + B^{(\sigma)}   y(t - \tau), \\ 
\nonumber y(t)&=&C^{(\sigma)} f(x(t)) + D^{\sigma} y(t-\tau)
\end{eqnarray}
and the corresponding family of subsystems
\begin{eqnarray}\label{eq:FamCoupled}
\dot{x} &=& A^{(s)} f(x(t)) + B^{(s)}   y(t - \tau), \\
\nonumber y(t)&=& C^{(s)}  
f(x(t)) + D^{(s)}  y(t-\tau), \quad s=1,\ldots,N.
\end{eqnarray}
It is not too difficult to adapt the argument of Theorem \ref{thm:Coupled1} to derive the following result. 
\begin{theorem}\label{thm:SwCoupled}
Consider the switched coupled system (\ref{eq:SwCoupled}).  Assume that $A^{(s)}$ is a Metzler matrix,  
and that $B^{(s)}, C^{(s)}, D^{(s)}$ are nonnegative matrices for $s \in \{1, \ldots, N\}$.  Further assume that there exists
a vector $v \gg 0$ such that 
\begin{equation}\label{eq:Dsw} 
{D^{(s)}}^T v \ll v \quad \mbox{ for } s = 1, \ldots, N.
\end{equation}  
If there exists a vector $\nu \gg 0$  satisfying
\begin{equation}\label{eq:swnu1}
(A^{(s)} +B^{(r)}(I-D^{(r)})^{-1}C^{(s)})^T\nu\ll 0,  \quad
                                      s,r=1,\ldots,N, 
\end{equation}
then we can choose positive real numbers $\mu_1, \ldots , \mu_n$ such 
that (\ref{eq:LKCoupled}) is a common L-K functional for the family (\ref{eq:FamCoupled}) and the system (\ref{eq:SwCoupled}) is absolutely stable. 
\end{theorem}

\textbf{Proof:} 
For $\varepsilon > 0$, define 
$$\mu = \textrm{max}_{r} (I-D^{(r)})^{-T} {B^{(r)}}^T \nu + \varepsilon v.$$
It can be immediately verified by direct calculation that 
\begin{equation}\label{eq:SwCoupB}
{B^{(r)}}^T \nu + (D^{(r)} - I)^T \mu \leq \varepsilon (D^{(r)} - I)^T v \ll 0
\end{equation}
for $r = 1, \ldots, N$.
On the other hand, it follows from our choice of $\nu$ that for $s, r = 1, \ldots, N$,
\begin{equation}\label{eq:SwCoupA}
{A^{(s)}}^T \nu + {C^{(s)}}^T (I - D^{(r)})^{-T} {B^{(r)}}^T \nu \ll 0. 
\end{equation}
Taking the elementwise maximum over $r$, it now follows that we can choose $\varepsilon > 0$ sufficiently small to ensure that
\begin{equation}\label{eq:SwCoupA2}
{A^{(s)}}^T \nu + {C^{(s)}}^T \mu \leq -\beta e
\end{equation}
for some $\beta >0$.  

If we differentiate the functional (\ref{eq:LKCoupled}) with respect to the $s$th subsystem from the family (\ref{eq:FamCoupled}), 
we obtain                                          
\begin{equation}\label{eq:DiffSwCoup}
\partial^+ V \leq    |f(x(t))|^T \left({A^{(s)}}^T \nu + {C^{(s)}}^T \mu \right) 
\end{equation}
$$
+ |y(t-\tau)|^T \left( {B^{(s)}}^T \nu + (D^{(s)}-I)^T \mu \right).
$$
It follows from (\ref{eq:SwCoupB}) and (\ref{eq:SwCoupA2}) that 
$$ \partial^+ V \leq - \beta \sum_{i=1}^n |f_i(x_i)|.$$
The result now follows immediately.

In the case where the switching only occurs in the state $x(t)$, so that both $B$ and $D$ are fixed, the conditions of the previous result can be relaxed somewhat.
\begin{corollary}\label{cor:Coupled}
 Consider the system (\ref{eq:SwCoupled}) and assume that $B^{(s)} = B$, $D^{(s)} = D$, where $B, D$ are nonnnegative and $D$ is Schur-Cohn for $s = 1, \ldots, N$.  Furthermore, suppose that $A^{(s)}$ is a Metzler matrix, and that  
 $C^{(s)}$ is nonnegative.  If there exists a positive vector $\nu$ satisfying
$$
(A^{(s)} +B(I-D)^{-1}C^{(s)})^T\nu\ll 0,  \quad
                                      s=1,\ldots,N,
$$
the system (\ref{eq:SwCoupled}) is absolutely stable,
and there exists a common 
L-K functional  
of the form (\ref{eq:LKCoupled}) for the family (\ref{eq:FamCoupled}).
\end{corollary}

\textbf{Remark:} As in Section \ref{sec:swit}, the results of this section can be extended to coupled systems with several delays.

\section{Neutral Type Systems}
\label{sec:neut}
Neutral systems, in which delays appear in the derivative of the state vector as well as in state itself arise in a number of applications and have been widely studied \cite{Hale}.  We next present a preliminary result on the stability of neutral systems using the techniques of the previous sections.  

Consider the switched linear neutral delay system 
\begin{equation}\label{eq:NeutralSys}
\dot{x}(t)-D\dot x(t-\tau) = A^{(\sigma)}x(t) + G^{(\sigma)}x(t - \tau) 
\end{equation}

and  the corresponding family of subsystems
\begin{equation}\label{eq:NeutralFam}
\dot{x}(t)-D\dot x(t-\tau) = A^{(s)}x(t) + G^{(s)}x(t - \tau), \ \
s=1,\ldots,N. 
\end{equation}

Here  $x\in \mathbb{R}^n$, 
$A^{(s)},G^{(s)},D$ are constant matrices in $\mathbb{R}^{n \times n}$ and $\tau$ is a constant nonnegative 
delay.  In the stability analysis of neutral systems, it is usual to assume that the operator on the derivative is stable.  In our context, this amounts to assuming that the matrix $D$ is Schur-Cohn stable.  We make this assumption from here on.

It is known, see for example \cite{Pepe2},
that the systems in the family (\ref{eq:NeutralFam}) can be transformed  into the family of 
coupled delay differential 
and continuous time difference systems
$$
\dot{y}(t) = A^{(s)} y(t) + B^{(s)}x(t - \tau), \qquad x(t)=y(t) + Dx(t-\tau),
$$
where $B^{(s)}=A^{(s)}D+G^{(s)}$.
 
With this in mind, it is possible to analyse the stability of the system (\ref{eq:NeutralSys}) using the results of the previous section.   
However, we shall use results presented in \cite{Hale} to provide a sufficient condition for the systems (\ref{eq:NeutralFam}) to have a common L-K functional.  
Note that, in keeping with a standard approach to neutral systems, the systems in the family (\ref{eq:NeutralFam}) can be rewritten as follows
\begin{equation}\label{eq:NeutralFam2}
\dot{x}(t)-D\dot x(t-\tau) = A^{(s)}(x(t)-D x(t-\tau)) + B^{(s)}x(t - \tau), 
\end{equation}        
for $s=1,\ldots,N$.  Here $B^{(s)}=A^{(s)}D+G^{(s)}$.

In contrast to previous sections, we will not assume that the above systems are positive.  In fact, characterising a positive neutral system is not straightforward, and even if we assume that the matrices $A^{(s)}$ are Metzler, and $D, G^{(s)}$ are nonnegative, then it does not necessarily follow that $B^{(s)}$ is nonnegative.  

Instead, we use properties of Metzler and nonnegative matrices in an indirect way.  We construct the auxiliary matrices  $\tilde A^{(s)}, \tilde B^{(s)}, \tilde D$
whose entries are given by
$$
\tilde a^{(s)}_{ii}=a^{(s)}_{ii}, \ \  \tilde a^{(s)}_{ij}=|a^{(s)}_{ij}| \ \ 
{\rm for}\ \    i\neq j,
$$
$$
\ \  \   \tilde b^{(s)}_{ij}=|b^{(s)}_{ij}|, \ \  \tilde d_{ij}=|d_{ij}|.
$$
With this notation, we have the following result.  

\begin{theorem}\label{thm:NeutralMain}
Consider the system (\ref{eq:NeutralSys}).  Suppose that $ \tilde D$ is a Schur-Cohn stable matrix, and 
that there exists a positive vector $\nu=(\nu_1,\ldots,\nu_n)^T$ 
satisfying
$$
(\tilde A^{(s)} +\tilde B^{(r)}(I-\tilde D)^{-1})^T\nu\ll 0,  \quad
                                      s,r=1,\ldots,N.
$$                                       
Then the switched system (\ref{eq:NeutralSys}) is asymptotically stable for arbitrary admissible switching law 
and for any $\tau\geq 0$.
\end{theorem}

\textbf{Proof:} 
Choose some vector $v \gg 0$ such that $(\tilde D - I) v \ll 0$ and for $\varepsilon > 0$, define 
$$\mu = \textrm{max}_s \left(\tilde B^{(s)}(I- \tilde D)^{-1} \right)^T \nu + \varepsilon v.$$
As in the proofs in the last two sections, it can be readily checked that 
$$(\tilde B^{(s)})^T \nu + (\tilde D - I)^T \mu \ll 0$$
for all $s$ and that by choosing $\varepsilon > 0$ sufficiently small, we can ensure that 
$$(\tilde A^{(s)})^T \nu + \mu \leq - \beta e$$
for some $\beta > 0$.

Now consider a L-K  functional in the form
$$
V=\sum_{i=1}^n \nu_i \left|x_i(t)-\sum_{k=1}^n d_{ik}x_k(t-\tau)\right|
$$
$$
+ 
\sum_{j=1}^n \mu_j
\int_{t-\tau}^t |x_j(z)|dz.
$$
The Dini upper-right derivative along the $s$th subsystem of (\ref{eq:NeutralSys}) of this functional satisfies
$$
\partial^+ V\leq \sum_{j=1}^n   \left|x_j(t)-\sum_{k=1}^n d_{jk}x_k(t-\tau)\right|
\left( \sum_{i=1}^n  \nu_i \tilde a^{(s)}_{ij} +   \mu_j \right)
$$
$$
+ \sum_{j=1}^n   |x_j(t-\tau)| \left(\sum_{i=1}^n  \nu_i \tilde b^{(s)}_{ij} 
+\sum_{i=1}^n  \mu_i \tilde d_{ij}-\mu_j\right).
$$
It follows from the first part of the proof that
$$
\partial^+ V\leq -\beta \sum_{j=1}^n   \left|x_j(t)-\sum_{k=1}^n d_{jk}x_k(t-\tau)\right|.
$$
The result now follows from Theorem 9.8.4 of \cite{Hale}.

\section{Switched Nonlinear Difference Systems with Time-Delay}
\label{sec:disc}
In this section, we briefly note that results on \emph{linear} functionals similar to those given in Section \ref{sec:swit} can also be derived for discrete-time systems.  As noted in \cite{KasBha}, the systems we consider here 
arise in applications such as digital filtering. 

As in Section {\ref{sec:swit}, we first consider the single delay case in the 
interest of notational simplicity and clarity.  

\subsection{A Single Delay}

Consider a switched system of the form
\begin{equation}\label{eq:discsw}
x(k+1)=A^{(\sigma)}f(x(k))+B^{(\sigma)}f(x(k-m)), 
\end{equation}
and the corresponding family of subsystems
\begin{equation}\label{eq:discfam}
x(k+1)=A^{(s)}f(x(k))+B^{(s)}f(x(k-m)), 
\end{equation}
for $s=1,\ldots,N$. Here, $A^{(s)},B^{(s)}$ are nonnegative matrices and $m$ is a nonnegative 
integer delay.
The nonlinearities $f(x)=(f_1(x_1),\ldots,f_n(x_n))^T$ are continuous for $x\in \mathbb R^n$ and
satisfy (\ref{eq:nonlin}).  In addition, in this section we will assume that
\begin{equation}\label{eq:discnonlin}
|f_i(x_i)|\leq |x_i|, \quad i=1,\ldots,n.
\end{equation}

We denote by $x_k$ the state of the delayed system given by $x_k = (x(k), x(k-1), \ldots, x(k-m))^T$.
\begin{theorem}\label{thm:DiscSwit1}
Let $A^{(s)}$, $B^{(s)}$ be nonnegative matrices for $s = 1, \ldots , N$.  Assume that there exists a 
vector $\nu \gg 0$ satisfying the
inequalities
\begin{equation}\label{eq:discswit1}
(A^{(s)}+B^{(r)}-I)^T \nu \ll 0, \qquad s,r=1,\ldots ,N.
\end{equation}
Then there exist positive numbers $\mu_1, \ldots, \mu_n$ such that 
\begin{equation}\label{eq:LKdisc}
V=\sum_{i=1}^n \nu_i |x_i(k)|+\sum_{i=1}^n \mu_i \sum_{l=1}^m 
|f_i(x_i(k-l))|
\end{equation}
defines a common L--K functional for the systems 
 (\ref{eq:discfam}).
\end{theorem}
                    
\textbf{Proof:}  
As our system is discrete-time, we consider the difference 
$\Delta V = V(x_{k+1}) - V(x_k)$ of the functional (\ref{eq:LKdisc}) with 
respect to the system (\ref{eq:discfam}) for some $s \in \{1, \ldots , N\}$.  Using (\ref{eq:discfam}) and (\ref{eq:discnonlin}), we can show by direct calculation that 
\begin{equation}
\label{eq:DeltaSw}
\Delta V\leq \sum_{j=1}^n |f_j(x_j(k))|\left(
\sum_{i=1}^n \nu_i a^{(s)}_{ij}-\nu_j+\mu_j\right) 
\end{equation}
$$
\sum_{j=1}^n |f_j(x_j(k-m))|   \left(
\sum_{i=1}^n \nu_i b^{(s)}_{ij}-\mu_j\right).
$$
It is clear from (\ref{eq:DeltaSw}) that if we can find a positive vector $\mu \gg 0$ such that 
\begin{equation}\label{eq:disccommon1}
(A^{(s)}-I)^T \nu  +\mu \ll 0, \quad {B^{(s)}}^T \nu-\mu \ll 0, \quad s=1,\ldots,N, 
\end{equation} 
then (\ref{eq:LKdisc}) will define a common L-K functional for the family of systems (\ref{eq:discfam}).

To see that such a vector $\mu$ must exist, note that as there are only 
finitely many inequalities in (\ref{eq:discswit1}), there exists some 
positive vector $w \gg 0$ with 
$$(A^{(s)} - I)^T \nu + (B^{(r)})^T \nu \ll -w$$
for $s, r$ in $\{1, \ldots, N\}$.  Define $d = \textrm{max}_{r}(B^{(r)})^T \nu$ 
as the componentwise maximum of the vectors $(B^{(r)})^T\nu$. It is now easy 
to verify that $\mu = d + \frac{w}{2}$ satisfies (\ref{eq:disccommon1}) and hence for this choice of $\mu$, (\ref{eq:LKdisc}) defines a L-K functional for (\ref{eq:discsw}) as claimed. 

\subsection{Several Delays}
Finally, for this section, consider the family of subsystems
\begin{equation}\label{eq:discmult}
x(k+1)=A^{(s)}f(x(k))+B^{(s)}_1f(x(k-1))
\end{equation}
$$
+\ldots+
B^{(s)}_mf(x(k-m)),
$$
for $s=1,\ldots,N$.  Here $A^{(s)},B^{(s)}_1,\ldots,B^{(s)}_m$ are nonnegative matrices for $1 \leq s \leq N$, 
and the nonlinearities $f(x)=(f_1(x_1),\ldots,f_n(x_n))^T$ satisfy the 
conditions (\ref{eq:nonlin}), (\ref{eq:discnonlin}).  

We note that it is relatively straightforward to adapt the argument used in the proof of Theorem \ref{thm:DiscSwit1} to obtain the following result for systems with several delays.
\begin{theorem}\label{thm:discmult}
Assume there exists a positive vector $\nu\gg 0$ satisfying the
inequalities
$$
(A^{(s)}+B^{(r_1)}_1+\ldots+B^{(r_m)}_m-I)^T \nu \ll 0, 
$$
for $ 
s,r_1,\ldots r_m =1,\ldots,N$.
Then there exist positive constants $\mu_{lj}$, such that 
\begin{eqnarray*}
V&=&\sum_{i=1}^n \nu_i |x_i(k)| \\
&+&\sum_{j=1}^n \sum_{l=1}^m \mu_{lj}  
\left(|f_j(x_j(k-1))|+\ldots+|f_j(x_j(k-l))|\right) 
\end{eqnarray*}
defines a common L-K functional for the family (\ref{eq:discmult}).
\end{theorem}

\section{Conclusions}
\label{sec:conc}   
We have presented a set of  results concerning Lyapunov-Krasovski functionals and absolute stability for various classes of nonlinear, switched systems with time-delay.  Specifically, we have described sufficient conditions for absolute stability of the system class considered in \cite{SunWang} that relax the requirements of this previous paper.  We have extended this analysis to systems described by coupled differential-difference equations and neutral systems.  We have also briefly noted that corresponding results can be obtained for discrete-time systems.

\end{document}